\documentclass[twoside,draft,reqno,12pt]{amsart}
\usepackage[english]{babel}
\usepackage{amsmath}
\usepackage{amssymb}
\usepackage{enumerate}

\usepackage{amsfonts}
\usepackage{graphicx}

\usepackage[ citecolor=black, linkcolor=black,
colorlinks, hypertexnames=false, plainpages=false]{hyperref}

\usepackage{color}

\newtheorem{lemma}     {Lemma}[section]
\newtheorem{thm}   [lemma]{Theorem}
\newtheorem{teorema1}   [lemma]{Theorem}
\newtheorem{prop}       [lemma]{Proposition}
\newtheorem{coro}       [lemma]{Corollary}
\newtheorem{cong1}      [lemma]{Conjecture}
\newtheorem{remark1}    [lemma]{Remark}
\newtheorem{defin}      [lemma]{Definition}

\numberwithin{equation}{section}

     \newcommand{\nn}{\nonumber}

\newcommand{\dis}{\displaystyle}

\newcommand{\mmmintone}[1]{{\dis{\int\kern -.38cm
-}}_{\kern-.21cm\substack{#1}}\;\;}
\newcommand{\mmmintwo}[2]{{\dis{\int\kern -.43cm
-}}_{\kern-.21cm\substack{#1}}^{\substack{#2}}\;\;}
\newcommand{\submint}{{\scriptstyle{\int\kern -.66em -}}}
\newcommand{\submintone}[1]{{\scriptstyle{\int\kern -.66em
-}}_{\scriptscriptstyle{\kern-.21em\substack{#1}}}}
\newcommand{\fracmint}{{\textstyle{\int\kern -.88em -}}}
\newcommand{\fracmintone}[1]{{\textstyle{\int\kern -.88em
-}}_{\scriptscriptstyle{\kern-.21em\substack{#1}}}\;}

\newcommand{\eps}{\epsilon}

\newcommand{\E}{\mathbb E}

\newcommand{\nada}[1]{}

\begin{document}
\today

\vskip.5cm
\title[Escape strategies]
{First passage percolation and escape strategies}

\author{E. D. Andjel}
\address{ Enrique D. Andjel, LATP URA 225/CNRS,\newline Universit\'e
d'Aix-Marseille, \newline
\indent 39, rue Joliot Curie. 13453 Marseille cedex 13, France}
\email {Enrique.Andjel@cmi.univ-mrs.fr}

\author{M.E. Vares}
\address{ Maria Eulalia Vares, \newline DME, Instituto de Matem\'atica, Universidade Federal do Rio de Janeiro, \newline
Av. Athos da Silveira Ramos 149, CEP 21941-909 - Rio de Janeiro, RJ,
Brasil }\email{eulalia@im.ufrj.br}
\begin{abstract}

Consider first passage percolation on $\mathbb{Z}^d$ with passage
times given by i.i.d. random variables with common distribution $F$.
Let $t_\pi(u,v)$ be the time from $u$ to $v$ for a path $\pi$ and
$t(u,v)$ the minimal time among all paths from $u$ to $v$. We
ask whether or not there exist points $x,y \in \mathbb{Z}^d$ and a
semi-infinite path $\pi=(y_0=y,y_1,\dots)$ such that $t_\pi(y,
y_{n+1})<t(x,y_n)$ for all $n$. Necessary and sufficient conditions
on $F$ are given for this to occur. When the support of $F$ is unbounded, 
we also obtain results on the number of edges with large passage time used 
by geodesics. 

\end{abstract}

\maketitle


\section{Introduction}

\label{sec0}

This work can be motivated by the following {\it game}: two
individuals, called $\lambda$ and $\sigma$, move on $\mathbb{Z}^d$,
spending a random time $\tau(e)$ to cross each edge $e=\langle x,y
\rangle$ from $x$ to a nearest neighbor site $y$, or vice-versa.
These passage times are assumed to be i.i.d. non-negative random
variables (further assumptions will be made later). The two
individuals start from distinct positions, that we denote by
$x_\lambda$ and $x_\sigma$ respectively; $\lambda$ would like to
catch $\sigma$, who in turn wants to escape. A natural and simple
question is: knowing the passage times, is  $\sigma$ able to devise
a strategy that would be successful for his/her goal, independently
of what $\lambda$ does? In this case we shall say that $\sigma$ has
a {\it perfect} strategy. Thus, a perfect strategy is an infinite
sequence of moves for $\sigma$ in a way that (s)he will never be
caught by $\lambda$ regardless of what (s)he does. We are imagining
that the individuals are always sitting on a vertex of
$\mathbb{Z}^d$ and think of an edge as being a door that remains
closed unless someone {\it knocks} at it.  To {\it open} a door corresponding to edge $e$ an individual
must be at one of the endpoints of $e$ and knock at it. Once this is done, the door will open after a time interval of length $\tau(e)$,
 the individual will cross that edge
and the door will be closed immediately after. To prevent trivial
situations it is natural to assume the common distribution to be
{\it useful} in the sense of Definition \ref{useful} below; in
particular $\{\tau(e)=0\}$ does not percolate.  As we shall see the
situation changes significantly depending on the distribution of the
passage times being supported on a compact set or not. When $F$ has
unbounded support, the probability that a perfect strategy for
$\sigma$ exists is zero, independently of the starting positions. On
the other hand, in the bounded case and, if a clairvoyant $\sigma$
can chose the initial position depending on $x_\lambda$ and the
$\tau$ variables, then a perfect strategy can be implemented with
probability one. The result is precisely stated in Theorem
\ref{thm1.1} below. The proof in the bounded case is very simple.
The unbounded case involves Proposition \ref{prop3.1}, which is
the key result of the paper. From it we also derive Theorem \ref{thm2.1} 
which states that for $F$ with unbounded support and any arbitrary 
constant $M$, geodesics use asymptotically a positive proportion of 
edges whose passage times are larger than $M$.

\vskip 0.1cm We now outline the paper: in the next section we
introduce the basic notation and definitions, and state two
theorems. Theorem \ref{thm1.1} has to do with the above question,
and its easy part is proved in the same section. Theorem
\ref{thm2.1} adds information to the comparison of first passage
percolation models, treated by van den Berg and Kesten in \cite{bk}
(in a more general case) in terms of time constants. In section
\ref{sec3} we state and prove Proposition \ref{prop3.1} which is the
main technical result of the paper, from which
 part (i) of Theorem \ref{thm1.1} and Corollary \ref{co2.1} follow immediately.
We then conclude the proofs and discuss a related
problem.

\section{Preliminaries and results}

\label{sec2}

\noindent {\bf Notation and definitions.}

In this paper $\E=\E^d$ will denote the set of nearest  neighbor
(n.n.) edges in the cubic lattice $\mathbb{Z}^d$. The origin in
$\mathbb{Z}^d$ will be denoted by $\mathbf{0}$.
For $x,y \in \mathbb{Z}^d$, $\|x-y\|$ will denote the
$\ell_1$--distance, i.e. $\|x\|=\sum_{i=1}^d|x^i|$ for
$x=(x^1,\dots,x^d)\in \mathbb{Z}^d$. A finite path
$\pi=(e_1,\dots,e_k)$ is a sequence of adjacent edges (sharing a
vertex), i.e. $e_i=<x_{i-1}, x_{i}>$ for each $i=1,\dots,k$. 
In this case we say
that $\pi$ goes from $x_0$ to $x_{k}$. For the context of this
paper, it suffices to consider self-avoiding paths, i.e. when the
$x_i, i=0,\dots,k$ are all distinct, and we always assume this
without further comment. Sometimes we identify a path with the
sequence of its visited vertices, writing $\pi=(x_0,\dots,x_k)$.

The basic random object consists of a family $\{\tau(e)\colon e \in \E\}$ of i.i.d. non-negative random variables
defined on a probability space $(\Omega, \mathcal F, P)$,  where $\tau(e)$ represents the passage time
of the edge $e$, interpreted as the time to traverse $e$. Their common distribution will be denoted by $F$.
The passage time $t(\pi)$ of a given path $\pi=(e_1,\dots,e_k)$ is simply given by the sum of the variables
$\tau(e_i)$ for $i=1,\dots,k$. We say that a given path $\tilde \pi$ from $x$ to $y$ is optimal (from $x$ to $y$)
if its travel time
is the shortest among all paths from $x$ to $y$:
\begin{equation}
\label{time}
 t(\tilde\pi)= \inf \{t(\pi)\colon \pi \text{ is a path from } x \text{ to } y\}=:t(x,y).
\end{equation}

Any such optimal path is also called a {\it geodesic} (from $x$ to
$y$). An infinite path $\tilde \pi=(e_1,e_2,\dots)$ starting at $x$
is said to be a semi-infinite geodesic if for any $n$ the finite
path $(e_1,\dots,e_n)$ is a geodesic from $x$ to its endpoint. It is
easy to see that semi-infinite geodesics starting from any given
point always exist. We also see easily that when $F$ is continuous,
there is a.s. a unique optimal path from $x$ to $y$ for any two distinct
vertices $x$ and $y$. Here, some assumptions on $F$ will be needed,
and as in \cite{bk} we set the following:

\begin{defin}
\label{useful}
A distribution $F$ with support in $[0,+\infty)$ is called {\em useful} if the following holds:

\begin{eqnarray}
\label{eq-useful}
F(r)< p_c \text{ when } r=0,\\\nn
F(r)< \overrightarrow{p_c}  \text{ when } r>0,\\\nn
\end{eqnarray}

where $p_c$ ($\overrightarrow{p_c}$) denotes the critical
probability for the Bernoulli (oriented, resp.) bond percolation
model on $\mathbb{Z}^d$, and $r$ stands for the minimum of the
support of $F$, hereby denoted by $supp(F)$.
\end{defin}

 \begin{thm}
                \label{thm1.1}

\noindent Let $F$ be useful in the sense of Definition \ref{useful}.
\newline
\noindent (i)\, If $F$ has unbounded support, then for any $x_\lambda, x_\sigma$
\begin{equation}
\label{eq1.2} P(\sigma \text{  has perfect strategy })=0.
\end{equation}
\noindent (ii)\, Assume $F$ to be supported in $[0,M]$ for some
finite $M$. Let $\tilde \pi$ be a semi-infinite geodesic from
$x_\lambda$. If the event
\begin{equation}
\label{eq1.1}
[M + t(x_\sigma,x) <t(x_\lambda, x) \text{ for some }x \in \tilde \pi]
\end{equation}
occurs, then  $\sigma$ has a perfect strategy. In particular, given
$x_\lambda$, with probability one there exist (infinitely many)
random initial positions $x_\sigma$ from where $\sigma$ has a
perfect strategy.

                \end{thm}

\vskip.5cm

\noindent {\bf Proof.} Part {\it (i)} will be proven in Section
\ref{sec3} as a corollary of our Proposition \ref{prop3.1}.  We now
prove only the easy part {\it(ii)}. Indeed, under the situation
described in \eqref{eq1.1}, it follows at once that a perfect
strategy for $\sigma$ consists in taking any $x \in \tilde\pi$ for
which $M+ t(x_\sigma,x)<t(x_\lambda,x)$, moving to $x$ by the
geodesics from $x_\sigma$ to $x$ and then following the infinite
branch of $\tilde\pi$ that starts in $x$. On the other hand,  if $F$
is useful it follows at once from the definitions that there exists
$\delta>0$ so that $F(\delta)<p_c$, which implies $t(x_\lambda,x)
\to \infty$ as $x \to \infty$ along $\tilde \pi$, and the inequality
in \eqref{eq1.1} becomes trivial for $x_\sigma \in \tilde \pi$ with
$t(x_\lambda,x_\sigma) >M$.

\medskip
We now state our second theorem.

\begin{thm}
\label{thm2.1} Let $F$ be a useful distribution on $[0,\infty)$ with
unbounded support. Then, for each $M$ positive there exists
$\eps=\eps(M)>0$ and $\alpha=\alpha(M)>0$ so that for all $n\ge 1$ and all $x$ with
$\|x\|=n$, we have
\begin{equation}
 \label{eq-teo1}
P\left(\exists\,  \text{ geodesic }\pi \text{ from } \mathbf{0} \text{ to } x
\text{ such that } \sum_{e \in \pi}{\bf 1}_{\{\tau(e) > M\}}  \leq \alpha n\right)\leq
e^{-\eps n}.
\end{equation}
\end{thm}


\section{Basic proposition. Proofs.}
\label{sec3}

\begin{prop}
\label{prop3.1}
Let $F$ be a useful distribution on $[0,\infty)$
with unbounded support. For each $M>0$ let
\begin{equation}
\label{3-1} \bar t_{(M)}(\mathbf{0},x)=\inf\{t(\pi)\colon \pi \text{
is a path from $0$ to $x$ with }\tau(e)\le M \; \text{ for all $e$
in $\pi$}\},
\end{equation}
with the understanding that $\inf \emptyset=+\infty$. Then, for each
$M$ positive there exists $\eps=\eps(M)>0$ and $n_0=n_0(M)$ so that
for all $n\ge n_0$ and all $x$ such that $\|x\|=n$, we have
\begin{equation}
\label{basico} P\left(  M+ t(\mathbf{0},x) <\bar
t_{(M)}(\mathbf{0},x)\right)\ge 1- e^{-\eps n}.
\end{equation}
\end{prop}

\vskip0.5cm
The proof of the above proposition  uses arguments from \cite{bk}
that are recalled below.
Before moving to it, we state  the following  immediate consequence:

\begin{coro}
\label{co2.1} Let $F$ be a useful distribution on $[0,\infty)$ with
unbounded support. Then, for each $M$ positive there exists
$\eps=\eps(M)>0$ so that for all $n\ge 1$ and all $x$ with
$\|x\|=n$, we have
\begin{equation}
 \label{eq-teo11}
P(\exists\,  \text{ geodesic }\pi \text{ from } \mathbf{0} \text{ to } x
\text{ such that } \tau(e) \le M \text { for all } e \in \pi) \le
e^{-\eps n}.
\end{equation}
\end{coro}

\vskip0.5cm
\noindent {\bf Notation.} For $N \in \mathbb N$ and $l=(l^1,\dots,l^d)\in \mathbb{Z}^d$ consider the following
partition of $\mathbb{Z}^d$ by hypercubes, as in \cite{bk} called $N$-cubes.
\begin{equation}
\label{cube} S_l(N)=\{x \in \mathbb{Z}^d\colon Nl^i \le x^i < Nl^i+
N, \forall i\}.
\end{equation}
The cubes are naturally indexed by $l$, and this indexing is also
used to define the distance between two $N$-cubes. If $C\subset
\mathbb{Z}^d$, we use $\mathcal{F}(C)$ to denote the $\sigma$-field
generated by the variables $\tau(e)$ corresponding to edges $e$ that
have both endpoints in the set $C$. \vskip 0.2cm The following
collections of boxes $T_l(N)$, $B^{+j}_l(N)$ and $B^{-j}_l(N)$ will
also be useful in the proofs: for $N \in \mathbb N,\, l \in
\mathbb{Z}^d$,
\begin{eqnarray}
 \label{large-c}
T_l(N)=\{x \in \mathbb{Z}^d\colon Nl^i -N \le x^i \le Nl^i+2N,
\forall i\}.\\\nn B^{\pm j}_l(N)=T_{l}(N) \cap T_{l\pm 2{\bf
e}_j}(N), \quad j=1,\dots, d,
\end{eqnarray}
where ${\bf e}_j, j=1,\dots, d$ denote the canonical unitary vectors.
\vskip 0.2cm
We first recall Lemma (5.2) from \cite{bk} (see also \cite{gk}) which follows from a Peierls argument:

\begin{lemma}
\label{lem1}
 If the
cubes $S_l(N)$ are colored black or white in a random fashion which
is (i) translation invariant; (ii) finite range (i.e. the color of
$S_l(N)$ is $\mathcal{F}(\cup_{l'}(S(l',N)\colon \|l'-l\| \le
c_0)$--measurable for a suitable constant $c_0$) and moreover,
$P(S_0(N) \text{ is black} ) \to 1$ as $N\to \infty$, then for all
$N$ sufficiently large we can find positive numbers
$\epsilon=\epsilon(N)$ and $D=D(N)$ so that for each $u,v \in
\mathbb{Z}^d$ the probability that each path from $u$ to $v$ visits
at least $\epsilon \|u-v\|$ distinct black $N$-cubes is not smaller
than $1-e^{-D\|u-v\|}$.
\end{lemma}

We must also recall Lemma (5.5) from \cite{bk}, which says that if
the distribution $F$ of the time variables is useful, then there
exist positive numbers $\delta=\delta(F)$ and $D_0=D_0(F)$ such that
\begin{equation}
 \label{1-1}
P(t(u,v) \le (r+\delta) \|u-v\|) \le e^{-D_0\|u-v\|},
\end{equation}
for all $u,v \in \mathbb{Z}^d$, where $r$ is as in Definition
\ref{useful}.

Although it is assumed troughout \cite{bk} that $F$ has finite first moment, this requirement is not used in the
proof of \eqref{1-1}.
\vskip 0.5cm Of course, for the proof of the proposition  it
suffices to consider $M>0$ large and such that $P(\tau(e) \in
(M,M+1])>0$. In the proof we shall also consider optimal paths for
the passage times

\begin{equation*}
\bar\tau(e)=\begin{cases} \tau (e)  &\text{ if  }\tau(e) \le M,\\
+\infty&\text{ otherwise.}
\end{cases}
\end{equation*}
\medskip

\noindent {\bf Black cubes.} We take $\delta=\delta(F)$ and
$D_0=D_0(F)$ so that \eqref{1-1} holds. Let $M \in (0, +\infty)$. We
now say that the $N$-cube $S_l(N)$ is {\it black} if for any path
$\pi$ lying entirely in $T_l(N)$ with endpoints $u,v$ such that
$\|u-v\| \ge N/4$ and using only edges with passage times less than
or equal to $M$, we do have $t(\pi) \ge (r+\delta) \|u-v\|$. The
$N$-cubes $S_l(N), S_{l'}(N)$ are said to be separated if $T_l(N)
\cap T_{l'}(N) =\emptyset$.

From \eqref{1-1} we see that Lemma \ref{lem1} applies. In
particular, having fixed $\delta$ as above, for any $N$ sufficiently
large we can take  $D=D(N,F)>0$ and $\epsilon=\epsilon(N,F)>0$ in a way
that for all $n$ large enough:
\begin{equation}
 \label{1-2}
P(\exists \text{ path from } \mathbf{0} \text{ to } \Gamma_n \text{
that visits at most } [\epsilon n] \text{ separated black }
N-\text{cubes}\} \le e^{-D n},
\end{equation}
where $\Gamma_n=\{x \in \mathbb{Z}^d \colon \|x\|=n\}$ and $[\cdot]$
denotes the integer part.  (Of course, changing $D$ we may
assume \eqref{1-2} holds for all $n$.)

We shall now work on the complement of the event on the l.h.s. of
\eqref{1-2}. For the proof of Proposition \ref{prop3.1} we will try
to improve over the optimal paths for $\bar\tau$ from $0$ to some
$x$ in $\Gamma_n$  by examining the probability of successful
shortcuts in disjoint boxes $B^{\pm j}_l(N)$. The main point is the
control of the conditional probability of a successful shortcut.

\begin{defin}
\label{crossing} We say that a path $\pi$ crosses the box $B^{\pm
j}_l(N)$ if it crosses the box in the shortest direction and, except
for its endpoints, is entirely contained in the interior of
$B^{\pm j}_l(N)$.
\end{defin}

In the sequel we write $\pi_{[u,v]}$ to denote the stretch of $\pi$
starting at $u$ and ending at $v$.

\begin{defin}
\label{shorcutable} We say that a stretch  $\pi_{[u,v]}$ of $\pi$ is
{\it shortcutable} if it crosses one of the boxes $B^{\pm j}_l(N)$
corresponding to a black cube $S_l(N)$.
\end{defin}
%
%
%
%
%

\begin{defin}
 For $\rho>0$ and small, we say that a path from the
origin to $\Gamma_n$ satisfies property ${\mathcal P}_n(\rho)$ if it
contains at least $[\rho n]$ (integer part of $\rho n$) shortcutable
stretches which lie at distance at least $14N$ of each other.
\end{defin}

\begin{lemma}
\label{propriedade-rho} Let $N$ be large enough for \eqref{1-2} to
hold. There exist constants $\rho=\rho(N,F)>0$ and $D=D(N,F)>0$ such
that for all $n$ the probability that all paths from the origin to
$\Gamma_n$ satisfy condition ${\mathcal P}_n(\rho)$ is at least
$1-\exp(-Dn)$.
\end{lemma}

\noindent {\bf Proof.} It is clear that if $\pi$ is a path
connecting a vertex $x$ in $S_l(N)$ to $y \notin T_l(N)$, it must
contain a path that crosses one of the $2d$  $N$--boxes $B^{\pm
j}_l(N)$ in the sense just defined. Hence the lemma follows at once
from \eqref{1-2}.\qed


\newpage

\noindent {\em Shortcuts}

Let $\pi$ be a path from $\mathbf{0}$ to a point in $\Gamma_n$ whose
edges have passage times less than or equal to $M$. Let $\pi^\prime$
be a shortcutable stretch of $\pi$ and call $B$ the $N$--box
(corresponding to a black cube) it crosses. Assuming $\pi$ to be
optimal for the $\bar \tau$ variables, we shall examine the
possibility of a successful shortcut $\tilde \pi$ for $\pi$ that
uses an edge with passage time larger than $M$. This would be a path
verifying the following conditions:

\begin{itemize}
\item $\tilde\pi$ and $\pi$ are edge disjoint;
\item the endpoints of $\tilde \pi$ coincide with those of a segment
$\pi^{\prime\prime}$ of $\pi$;
\item  $|\tilde \pi| \le c_dN$ where the positive constant $c_d$
depends only on the dimension;
\item $\tilde \pi$ is contained in the same $N$-box $B$ as
$\pi^\prime$.

\noindent We shall then say that a shortcut as above is {\em
successful} if $M+ t(\tilde \pi) < t(\pi^{\prime\prime})$.
\end{itemize}
\noindent Let us first assume for notational simplicity that $d=2$,
$B=B^1_l(N)$, which we write as $B=[a,a+N]\times [b,b+3N]$, and that
$\pi^\prime$ crosses $B$ from left to right. Writing
$\pi=(x_0,\dots,x_s)$,  let $v=x_j$ be the position in
$\{a+N\}\times [b+1,b+3N-1]$ where $\pi$ first reaches the rightmost
face of $B$ after entering $B$ and $u=x_i$ the position in $\{a\}
\times [b+1,b+3N-1]$ of the leftmost face of $B$ last visited before
getting to $v$, so that $i<j$ and $\pi^\prime$ is the segment of
$\pi$ that goes from $u$ to $v$, which we denote as $\pi_{[u,v]}$.
We choose $N=4K$ for some $K\in \mathbb N$. We may define as well
the vertex with lowest second coordinate and first coordinate in
$[a+K,a+3K]$ along $\pi_{[u,v]}$. If there are several such points,
let us take e.g. the leftmost one, call it $z=(z^1,z^2)$. We assume
that $z$ is on the leftmost half of $B$, i.e. $z^1\le a+N/2$
(the argument being analogous when $z$ is on the rightmost half of
$B$); we now define $\tilde \pi$ by starting from $z$ moving
downwards one step to $z^\prime= z-\mathbf{e}_2$ and then moving
horizontally to the right for at most $K$ steps or until we reach
any point $w$ in $\pi$, whatever comes earlier (note that we can
have $w=z^\prime$). In the first case, we then move vertically upwards until
reaching a vertex $w$ visited by $\pi$; this just defined
path from $z$ to $w$ is what we call $\tilde \pi$. Three cases have
to be analyzed:

(a) $w\in \pi_{[u,v]}$,

(b) $w$ is visited by $\pi$ before $u$,

(c)  $w$ is visited by $\pi$ after $v$.

In all of these three cases we define a new path substituting the stretch of $\pi$ between $w$ and $z$ by $\tilde \pi$

In case (a) the substituted part is the stretch $\pi_{[z,w]}$
contained in $\pi_{[u,v]}$. In this case $\|z-w\|\ge K$.

In case (b) the substituted part is the portion of $\pi$ going from $w$ to $z$.

In case (c) the substituted part is the portion of $\pi$ going from $z$ to $w$.

It is easy to check that in all three cases the substituted part of
$\pi$ contains a stretch of $\pi_{[u,v]}$ connecting two points at
distance at least $K=N/4$. If $\pi_{[u,v]}$ is shortcutable, then
the time of the substituted part is at least $\|z-w\|r +K\delta$.

The extension to higher dimension is simple and we always have
$\|z-w\| \le 3Nd=12Kd$.

Assuming that the stretch $\pi_{[u,v]}$ is shortcutable, a condition
which guarantees a successful shortcut is $M +\sum_{e_i \in \tilde
\pi}\tau (e_i)< \|z-w\|r+K\delta$. Note that the number of edges in
$\tilde \pi$ is at most $\|z-w\|+2$. For our application below
(proof of Proposition \ref{3-1}), we shall impose for one of the
edges, call it $e_1$, that $\tau(e_1) \in (M,M+1)$ and for the other
edges we impose passage times in the interval $[r,r+\delta^\prime)$
with $\delta^\prime=\delta/(24d)$. Hence, the shortcut is successful
if
\begin{equation*}
 2M+1+ (\|z-w\|+1)(r+\delta^\prime)< \|z-w\|r +K\delta,
\end{equation*}
which is implied by
\begin{equation} \label{sucesso} 2M+1+
r+\delta/(24d)+K\delta/2<K\delta,
\end{equation}
and this, in its turn, is satisfied for $K$ large enough (depending
on $M$, $\delta$
 and $r$).

\medskip
\noindent {\bf Proof of Proposition \ref{prop3.1}.}

Let $\delta$ be as in the definition of black cubes and take $M>0$
so that $P(\tau \in (M,M+1])>0$. We fix $N=4K$ large so that the
conclusion of Lemma \ref{propriedade-rho} holds and moreover
$2M+1+r+\delta/(24d)<K\delta/2$ as in the above construction.
 We may as well assume that the set on the right side of
\eqref{3-1} is not empty, and let $\Pi$ be a path where the minimum
is attained. (In case of non-uniqueness, the argument will apply to
any of the finitely many optimal paths, and any deterministic way to
list them will do the job.)

We now define random variables $U_1,V_1,U_2,V_2,\dots$ taking 
values in $\mathbb{Z}^d \cup \{\infty\}$. On the event
$\{\Pi=\pi\}$, $U_1,V_1 $ are such that $\pi_{[U_1,V_1]} $ is the
first shortcutable stretch of $\pi$. If no such stretch exists then
$U_1=V_1=\infty$; $U_2,V_2 $ are such that $\pi_{[U_2,V_2]} $ is the
first shortcutable stretch of $\pi$ after $V_1$ whose distance to
$\pi_{[U_1,V_1]} $  is at least $7N$. In general, $U_{i+1},V_{i+1}$
is such that $\pi_{[U_{i+1},V_{i+1}]}$ is the first shortcutable
stretch of $\Pi$ after $V_i$ whose distance to $\cup_{j=1}^i
\Pi_{[U_j,V_j]}$ is at least $7N$,  or $U_{i+1}=V_{i+1}=\infty$ if
no such stretch exists.

For a given  $n$ let $q=q(n)=[\rho n]$. Then, partition the
probability space in events as
$A(\pi,x_1,y_1,\dots,x_q,y_q)=\{\Pi=\pi,U_i=x_i,V_i=y_i: 1\leq i
\leq q\}$ and the event $G=\{U_q=+\infty\}$. For each of the
shortcutable stretches of $\pi$ there is a path $\tilde \pi_i$ as
defined above, with $z_i,w_i$ the corresponding vertices in that
construction.

Call $e_{i,1},\dots,e_{i,k_i}$ the edges of  $\tilde \pi_i$ and call
$e'_{i,1},\dots,e'_{i,{\ell_i}}$ the edges which have one endpoint
in $\tilde \pi_i \setminus \{w_i,z_i\}$ and whose other endpoint is
not in  $\tilde \pi_i$. We now define the event

\begin{eqnarray}
\label{eventoF}
&F_i(\pi,x_1,y_1,\dots,x_q,y_q)=A(\pi,x_1,y_1,\dots,x_q,y_q)\cap\\\nn
&\{\tau(e_{i,1})\in (M, M+1],
\tau(e_{i,2})<r+\delta^\prime,\dots,\tau(e_{i,k_i})
<r+\delta^\prime, \tau(e'_{i,1})>M, \dots,
\tau(e'_{i,{\ell_i}})>M\},
\end{eqnarray}
with $\delta^\prime=\delta/(24d)$ as before. Notice that
$k_i,\ell_i$ are uniformly (in $i$) bounded by a constant that
depends only on $K$ and $d$.

If the event $F_i(\pi,x_1,y_1,\dots,x_n,y_n)$  occurs then
substituting a part of $\pi$ as explained before we get a new path
$\pi'_i$, and $M+t(\pi'_i)<t(\pi)$.

Since by Lemma \ref{propriedade-rho} $P(G) \leq e^{-Dn}$, to
conclude the proof it suffices to show that
\begin{equation}
\label{atalho3} P(\cap_{i=1}^q F_i^c
(\pi,x_1,y_1,\dots,x_q,y_q)\vert A(\pi,x_1,y_1,\dots,x_q,y_q))\leq
(1-\varepsilon)^q \end{equation} for some $\varepsilon>0$
(independent of $q$).

The proof of \eqref{atalho3} will follow by a suitable application
of the following simple lemma.

\begin{lemma}
\label{lem2} Let $\Omega=\mathbb{R}^\Lambda$, where $\Lambda$ is a
finite or countable set, endowed with the usual product Borel
sigma-field $\sigma(\Lambda)$. Let $\Lambda_1$ be a (non-empty)
finite proper subset of $\Lambda$ and
$\Lambda_2=\Lambda\setminus\Lambda_1$. For $i=1,2$ let
$\Omega_i=\mathbb{R}^{\Lambda_i}$, so that $\Omega=\Omega_1
\times\Omega_2$ and $\sigma(\Lambda)=\sigma(\Lambda_1) \times
\sigma(\Lambda_2)$. Let $\mu_i$ be a Borel probability measure on
$(\Omega_i,\sigma(\Lambda_i))$, $i=1,2$ and $\mu=\mu_1\times\mu_2$
the product measure on $(\Omega,\sigma(\Lambda))$. If  $A \in
\sigma(\Lambda)$ and $\hat B\in \sigma(\Lambda_1)$ have the property
that  $x=(x_1,x_2) \in A$ and $y_1 \in \hat B$ imply $(y_1,x_2) \in
A$, then
$$\mu(B \cap A) \ge \mu(B)\mu(A),$$
where $B=\hat B\times \Omega_2$.


\end{lemma}

\noindent {\bf Proof.} The hypothesis on $A$ and $\hat B$ can be
written as

$$\mathbf{1}_{\hat B}(y_1)\mathbf{1}_A(x_1,x_2) \le \mathbf{1}_{\hat B}(y_1)\mathbf{1}_A(y_1,x_2)$$
for all $x_1,y_1 \in \Omega_1$ and all $x_2\in \Omega_2$. We compute
the iterated integral $\mu_1(dy_1)\mu_1(dx_1)\mu_2(dx_2)$ on both
sides. The left hand side yields, by Tonelli's theorem,
$\mu(B)\mu(A)$. On the right hand side we have
\begin{equation*}
\int_{\Omega_1}\mu_1(dy_1)\int_{\Omega_1}\mu_1(dx_1)\int_{\Omega_2}\mu_2(dx_2)
\mathbf{1}_{\hat B}(y_1)\mathbf{1}_A(y_1,x_2)
\end{equation*}
which again by Tonelli's theorem can be rewritten as
\begin{equation*}
\int_{\Omega}\mu(dy_1,dx_2)\mathbf{1}_{B\cap
A}(y_1,x_2)\int_{\Omega_1}\mu_1(dx_1)=\mu(B \cap A)
\end{equation*}
proving the lemma.\qed


\vskip .3cm

To conclude the proof of Proposition \ref{prop3.1}, let
$$A=A(\pi,x_1,y_1,\dots,x_q,y_q)=\{\Pi=\pi,U_i=x_i,V_i=y_i: 1\leq i
\leq q\}\ \mbox{and}$$
$$\hat B_i=\{\tau(e_{i,1})\in (M, M+1],
\tau(e_{i,2})<r+\delta^\prime,\dots,\tau(e_{i,k_i})
<r+\delta^\prime, \tau(e'_{i,1})>M, \dots,
\tau(e'_{i,{\ell_i}})>M\},$$ for $i=1,\dots,q$, so that $P(\hat
B_i)\ge \eta >0$ for all $i$. A few instants of reflection show that
the condition in the lemma is verified for the pair $A$ and $\hat
B_1$:  since $\Pi$ is optimal for the $\bar \tau$ variables and has
a finite time, it cannot cross any of the edges
$e'_{1,1},\dots,e'_{1,{\ell_1}}$; this prevents it from using the
advantageous edges $e_{1,2},\dots,e_{1,k_1}$, and therefore the
modified configuration remains in $A$. Call $F_i$ the event defined
in \eqref{eventoF}. Since $F_1=A \cap \hat B_1$, the lemma implies
that $P(F_1^c\mid A)\le 1-P(\hat B_1)$. Analogously, we can again
apply the lemma with $A$ replaced by $A \cap \cap_{j=1}^{i-1} F_j^c$
for $i=2,\dots,q$ and $\hat B_{i}$ to conclude that conditional
probability on the l.h.s. of \eqref{atalho3} is bounded from above
by $(1-\eta)^{[qn]}$. \qed

\medskip

\begin{remark1}\label{r1}
The argument used to prove Proposition \ref{prop3.1} also shows that
under the same conditions, for each $M$ positive there exist $\alpha>0$ and $\epsilon>0$ (both 
depending on $M$) so that
for all $n\ge 1$  and all $x \in \Gamma_n$, we have
\begin{equation}
\label{basico2} P\left(t(\mathbf{0},x)+ \alpha n <\bar
t_{(M)}(\mathbf{0},x)\right)\ge 1- e^{-\eps n}.
\end{equation}
\end{remark1}

Recall now  the time constant associated 
 to the passage time distribution $F$, given by the deterministic limit (in probability):
\begin{equation*}
\mu_F=\lim_{n \to \infty} \frac 1n t(\mathbf{0}, n\mathbf{e}_1).
\end{equation*}
The limit is also a.s. and in $L_1$ under conditions on the tail of
$F$, e.g. if $F$ has finite mean, and also holds along any fixed direction, 
not only the coordinate axis. (See \cite{R,CD,K,M}.)

When $F$ has exponentially decaying tail, the result in Remark \ref{r1}, taking e.g. 
$x=n\mathbf{e}_1$, can be seen as consequence of large deviation estimates on the variables $t(0,x)/n$. 
(See e.g.\cite{K}.)  
\medskip

\noindent {\bf Proof of part (i) of Theorem \ref{thm1.1}.}

It is clear that with probability one, no perfect strategy for
$\sigma$ can consist in remaining in a finite set for all times and
so it must reach the set $\{x\colon \|x-x_\sigma\|=n\}$ for any $n$. It
is obvious that on the event $\{t(x_\lambda,x_\sigma)< M\}$ any
perfect strategy can only use edges with passage time smaller than
$M$  and  cannot include finite paths between two points whose
passage time exceeds the minimal passage time between these points
by more than  $M$.  Thus
\begin{equation*}
P(\sigma \text{ has perfect strategy}, t(x_\lambda,x_\sigma) < M)
\le P(\exists x \colon \|x-x_\sigma\| =n, \bar
t_{(M)}(x_\sigma,x)\le M+ t(x_\sigma,x))
\end{equation*}

Given $\eta>0$, let $M$ be such that $P(t(x_\lambda,x_\sigma)\ge
M)\le \eta$. Given such $M$ we take $n$ and $\epsilon$ so that 
\eqref{basico} holds and  $c_dn^{d-1}e^{-\epsilon n} \le \eta$,
where $c_dn^{d-1}$ is an upper bound for the cardinality of
$\{x\colon \|x\|=n\}$. We then have
\begin{equation*}
P(\sigma \text{ has perfect strategy}) \le 2\eta.
\end{equation*}
\qed

\medskip

\noindent {\bf Proof  of Theorem \ref{thm2.1}.}
An $N$-cube $S_l(N)$ is now colored black if every geodesic from a point on its boundary to a point on the boundary of
$T_l(N)$ uses at least one edge whose passage time is larger than $M$. From Corollary \ref{co2.1} we see that Lemma \ref{lem1}
applies, and the result follows. \qed

\medskip

\begin{remark1}\label{r2}
The same argument used for Theorem \ref{thm2.1} yields the following result: Let $A$ be a Borel set to which $F$ attributes 
positive measure. Then,
\begin{equation*}
P\left(\exists\,  \text{ geodesic }\pi \text{ from } \mathbf{0} \text{ to } x
\text{ such that } \sum_{e \in \pi}{\bf 1}_{\{\tau(e) \in A \}}  \leq \alpha n\right) \leq
e^{-\eps n}.
\end{equation*}
\end{remark1}
This last result improves inequality (2.16) in \cite{bk}.

\bigskip

\noindent {\bf An open problem.} Suppose that we have a family of
independent Poisson processes $\{\mathcal{P}_e\colon e \in
\mathbb{E}^d\}$ of parameter $c>0$. Assume that the individuals
$\lambda$ and $\sigma$ may now move from a vertex $x$ to a n.n.
vertex $y$ at the jump times of $\mathcal{P}_e$ where $e=\langle
x,y\rangle$. In dimension one, as easily seen, the probability that
$\sigma$ has a perfect strategy is zero. What happens in higher
dimensions?

\bigskip

\noindent {\bf Acknowledgement.} We warmly thank a referee who showed us that Theorem \ref{thm2.1} could be derived from Corollary 
\ref{co2.1}. E.D.A. thanks the kind hospitality
of IMPA, Rio de Janeiro during the preparation of this paper. The
work of E.D.A. was partially supported by CNRS and by CAPES. M.E.V
is partially supported by CNPq grant 304217/2011-5 and project
474233/2012-0. Most of this work was done while M.E.V. was
researcher of CBPF.

\bibliographystyle{amsalpha}

\end{document}